\documentclass[a4paper,12pt]{scrartcl}

\usepackage[utf8]{inputenc}

\usepackage{amsmath,amsthm,graphicx,skmath}
\usepackage{bm}

\usepackage{graphicx,color}

\usepackage{bera}
\usepackage[charter]{mathdesign}

\usepackage{cite}
\usepackage{xspace,microtype}

\DeclareMathOperator{\sign}{sgn}
\newcommand{\diag}{\text{\textbf{diag}}}

\begin{document}


\title{Numerical investigation of Mach number consistent
  Roe solvers for the Euler equations of gas dynamics}
\author{Friedemann Kemm}

\maketitle



\paragraph*{Keywords:}
\label{sec:keywords}
Shock instability,
Carbuncle phenomenon,
Low Mach number effect,
Low-dissipation schemes,
Shock-capturing


\begin{abstract}
  While traditional approaches to prevent the carbuncle phenomenon in
  gas dynamics simulations increase the viscosity on entropy and shear
  waves near shocks, it was quite recently suggested to instead
  decrease the viscosity on the acoustic waves for low Mach
  numbers. The goal is to achieve what, in this paper, we call Mach
  number consistency: for all waves, the numerical viscosity decreases
  with the same order of the Mach number when the Mach number tends to
  zero. We take the simple approach that was used for the proof of
  concept together with the simple model for the increased numerical
  viscosity on linear waves and investigate the possibilities of
  combining both in an adaptive manner while locally maintaining
  Mach number consistency.
\end{abstract}








\section{Introduction}
\label{sec:introduction}

Since the seminal paper of Quirk~\cite{quirk}, an immense amount of
research has been conducted on the instability problem that appears in
the simulation of supersonic flow problems known as \emph{Carbuncle
  Phenomenon}. The name originates from the fact that in
strongly supersonic flows against an infinite cylinder simulated on a
body-fitted, structured mesh the middle part of the resulting bow
shock degenerates to a carbuncle-shaped structure. It was conjectured
already by Quirk~\cite{quirk} that this phenomenon is closely related
to other instabilities such as the so-called
\emph{odd-even-decoupling} encountered in straight shocks aligned with
the grid.
The research in this area was twofold. On the one hand, the stability
of discrete shock profiles was investigated in one as well as in
several space
dimensions~\cite{carb15,moschetta-vorticity,aroracarb,elling_carbuncles_2006,elling_carbuncle_2009,bultelle98,chauvat05,gressier00,helzbale,jinliu,michael-carbuncle,pandolfi,park-kwon,ramalho2010possible,robinet00,roe_carbuncles_2005,zhang-zhang,Kalkhoran200063,kitamura13}. On
the other hand, there was a lot of effort spent to find cures for the
failure of some schemes in numerical
calculations~\cite{huang-wu,kim03,pandolfi,park-kwon,sanders98,phong,loh_time-accurate_2007,zaide_shock_2011,zaide-diss,zaide_flux,tu2014evaluation,kitamura09,sangeeth-2018,chen-et-al-shock-stable,shock-stable-hllc,xie-num-instabilities,rodionov2d,rodionov3d,rodionov2019,mandal-robust-hll,garicano-mena-deconinck}.
Unfortunately, the failure is only found in schemes giving high
resolution of shear and entropy waves by treating them explicitly,
i.\,e.\ with so called \emph{complete Riemann solvers} (as opposed to
so called \emph{incomplete Riemann solvers}), a category that includes,
for example, the Godunov, Roe, Osher, HLLC and HLLEM
schemes~\cite{toro,hllem}.  Complete Riemann solvers are preferable in
calculations involving complex wave structures as well as boundary
layers.

It was found that even in one space dimension there are some
instabilities of discrete shock profiles: slowly moving shocks produce
small post-shock oscillations~\cite{quirk,aroracarb,jinliu}.
\index{post-shock oscillations} But also \index{shock!steady} in the
case of a steady shock, instabilities can be found depending on the
value of the adiabatic coefficient \(\gamma\) as was shown by Bultelle
et~al.~\cite{bultelle98}. These and their relationship to two-dimensional instabilities of discrete shock profiles were our main
focus in our earlier study on the sources of the carbuncle~\cite{carb15}. 

Another one-dimensional issue that was not yet considered in
connection with the carbuncle, an issue that is also not studied
in~\cite{carb15}, is the wrong amplitude of pressure fluctuations
normal to the flow direction that might arise from inconsistent
numerical viscosities. This effect was brought to attention quite
recently by Fleischmann et~al.~\cite{fleischmann}, who discuss the
influence of what we will call \emph{Mach number consistency}, namely
requiring the viscosity on each wave to be of the same
order~\(\mathcal O(M^\alpha) (M\to 0)\), on the stability of discrete
shock profiles.  While they force Mach number consistency in a strict
manner, i.\,e.\ \(\alpha = 2\), we loosen the definition of the term
by only requiring the viscosity to be of the same order for all
waves. In this sense, the traditional carbuncle cures are also Mach
number consistent in the vicinity of discrete shock profiles by
ensuring order~\(\mathcal O(M)\) for all viscosities.

In this study, we investigate the possibility of combining the two
approaches~\(\mathcal O(M^2)\) and~\(\mathcal O(M)\) via an adaptive
blending which would also allow~\(\mathcal O(M^{2-\beta})\) for some
parameter~\(\beta \in [0,1]\). As a starting point, we employ the
simple modified Roe schemes proposed by Fleischmann
et~al.~\cite{fleischmann}, which, although not leading to a production
ready scheme, still allow for valuable insights into the problem and
useful directions for further research.
%
In Section~\ref{sec:mach-numb-cons}, we discuss Mach consistent Roe
solvers, wherein the basic properties of the Roe scheme are reviewed
as well. Once a desired adaptive blending with numerical viscosity of
order~\(\mathcal O(M^{2-\beta})\) and a computationally cheaper
approximation by substituting the weighted geometric mean with a
weighted arithmetic mean is achieved, we perform some numerical tests,
cf.\ Section~\ref{sec:numer-invest}. These tests are chosen such that
the desired information on the overall behaviour of the schemes can be
extracted from their results. The ultimate goal as stated in
Section~\ref{sec:concl-poss-direct} is to construct a robust all Mach
number Roe solver that prevents the carbuncle, based on the insight
gained in this study.

\section{Mach number consistent Riemann solvers}
\label{sec:mach-numb-cons}

\subsection{The quasi linear form of the Euler equations}
\label{sec:quasi-linear-form}

The Euler equations for an inviscid gas flow are
\begin{align*}
  \rho_t + \nabla \cdot [\rho \vec v] & = 0\;, \\
  (\rho \vec{v})_t + \nabla \cdot [\rho \vec v \circ \vec v + p\vec I]  &  =
  0\;, \\ 
  E_t + \nabla \cdot [(E+p)\vec v] & = 0\;
\end{align*}
with the density \(\rho\), the velocity \(\vec v = (u,v,w)^T\) for the
3d-case and \(\vec v = (u,v)^T\) for the 2d-case, the pressure \(p\),
and the total energy \(E\). Density, velocity, and pressure are called
\emph{primitive variables} in contrast to the \emph{conserved
  variables} density, momentum, and total energy. Throughout this
study, we consider the case of an ideal gas. 
%
The system is rotationally\index{invariance!rotational} and Galilean
invariant. Thus, it is sufficient so consider the quasi linear form
\begin{equation}
  \label{eq:7}
  \vec q_x + \vec A^{(x)}\vec q_x + \vec A^{(y)}\vec q_y + \vec A^{(z)}\vec q_z\;,
\end{equation}
where~\(\vec A^{(x)}, \vec A^{(y)}, \vec A^{(z)}\) are the flux
Jacobians for the respective space directions, only for one space
dimension, i.\,e.~\(\pd{}y=\pd{}z =0\). Now, the eigenvalues
of~\(\vec A^{(x)}\) are the wave speeds in \(x\)-direction, and the
left and right eigenvectors determine the waves themselves.

Since we restrict our study to the two-dimensional case, we can drop
the equation for the momentum in \(z\)-direction and end up with the
wave speeds
\begin{equation}
  \label{eq:2}
  \lambda_1 = u + c\;,\qquad \lambda_2=\lambda_3=u\;,\qquad \lambda_4
  = u + c
\end{equation}
with the speed of sound
\begin{equation*}
  c = \sqrt{\frac{\gamma p}{\rho}} =
  \sqrt{(\gamma-1)\biggl(H-\,\frac{\vec v^2}{2}\biggr)}\;,  
\end{equation*}
where~\(\gamma\) is the ratio of specific heats and
\begin{equation}
  \label{eq:147}
  H = \frac{E+p}{\rho}
\end{equation}
the total enthalpy. The right eigenvectors of the flux Jacobian are
\begin{equation}
  \label{eqcc:26}
  \vec r_1 =
  \begin{pmatrix}
    1 \\ u-c \\ v \\ H -cu
  \end{pmatrix}\;,
  \vec r_2 =
  \begin{pmatrix}
    1 \\ u \\ v \\ \frac{1}{2}\,\vec v^2
  \end{pmatrix}\;, \quad
  \vec r_3 =
  \begin{pmatrix}
    0 \\ 0 \\ 1 \\ v
  \end{pmatrix}
  \vec r_1 =
  \begin{pmatrix}
    1 \\ u+c \\ v \\ H +cu
  \end{pmatrix}
\end{equation}
and the corresponding left eigenvectors 
\begin{equation}
  \label{eqcc:27}
  \begin{split}
    \vec l_1 & = \frac{1}{c\,(2H-\vec v^2)}\, \Bigl( Hu -
    \frac{(u-c)\,\vec v^2}{2},\ 
    \frac{\vec v^2}{2} - (H+cu),\ -cv,\ c\Bigr)\;,\\
    \vec l_2 & = \frac{1}{H-\frac{\vec v^2}{2}}\, \bigl( H-\vec v^2,\,
    u,\, v,\, -1\bigr)\;, \\ 
    \vec l_3 & = \bigl(-v,\, 0,\, 1,\, 0  \bigr)\;,\\
    \vec l_4 & = \frac{1}{c\,(2H-\vec v^2)}\, \Bigl( -Hu +
    \frac{(u+c)\,\vec v^2}{2},\ 
    -\frac{\vec v^2}{2} + (H-cu),\ -cv,\ c\Bigr)\;.
  \end{split}
\end{equation}
In the following, we denote the matrix with the left eigenvectors as
rows~\(\vec L\), the matrix with the right eigenvectors as
columns~\(\vec R\), and keep in mind that we have normalized the
eigenvectors such that~\(\vec L = \vec R^{-1}\).

\subsection{Roe's consistent local linearization}
\label{sec:roes-cons-local}

If we write the generic 1d-scheme
for a 1d conservation law
\begin{equation}
  \label{eq:140}
  \vec q_t + \vec f(\vec q)_x = \vec 0
\end{equation}
as
\begin{equation}
  \label{eq:5}
  \frac{\vec q_i^{n+1} - \vec q_i^n}{\Delta t} + \frac{\vec
    G_{i+1/2}^n - \vec G_{i-1/2}^n}{\Delta x} = \vec 0\;,
\end{equation}
with the generic numerical flux function
\begin{equation}
  \label{eq:4}
  \vec G (\vec q_r,\vec q_l) = \frac{1}{2} \bigl(f(\vec q_r) + f(\vec
  q_l)\bigr) - \frac{1}{2} \vec V (\vec q_r,\vec q_l) (\vec
  q_r-\vec q_l)\;, 
\end{equation}
it is obvious that the matrix~\(\vec V\) determines the numerical
viscosity. If we now apply this to a linear system
\begin{equation}
  \label{eq:33}
  \vec q_t + \vec A \vec q_x = \vec 0
\end{equation}
with constant system matrix~\(\vec A\) and employ standard upwind on
all waves, we end up with the viscosity matrix
\begin{equation}
  \label{eq:40}
  \vec V = \abs{\vec A}\;,
\end{equation}
with the convention
\begin{equation}
  \label{eq:41}
  \abs{\vec A} = \vec R\, \abs{\vec \Lambda}\, \vec L\;, \qquad \abs{\vec
    \Lambda} = \diag\,(\,\abs{\lambda_1},\dots,\abs{\lambda_m}\,)\;. 
\end{equation}
Thus, the absolute values of eigenvalues determine the numerical
viscosity.  This was the motivation for Roe~\cite{roe-orig} to look
out for consistent local linearizations, i.\,e.\ a
matrix~\(\tilde{\vec A} = \tilde{\vec A}(\vec q_l, \vec q_r)\) with the
following properties
\begin{gather}
  \label{eq:45}
  \vec f (\vec q_r) - \vec f (\vec q_l) = \tilde{\vec A} (\vec q_l, \vec
  q_r)\,(\vec q_r - \vec q_l)\;,  \\
  \label{eq:47}
  \tilde{\vec A} (\vec q_l, \vec q_r) \to \vec A(\vec q) \qquad
  \text{for}\quad (\vec q_l, \vec q_r) \to (\vec q, \vec q)\;,
   \\ 
  \label{eq:51}
  \tilde{\vec A} (\vec q_l, \vec q_r)\ \text{is diagonalizable for
    all}\ \vec q_l, \vec q_r\;.  
\end{gather}
A matrix~\(\tilde {\vec A} (\vec q_l, \vec q_r)\) that satisfies these
conditions is called a \emph{Roe matrix} or a \emph{consistent local
  linearization} for~\eqref{eq:140}. If there exists a single
state~\(\tilde{\vec q}=\tilde{\vec q}(\vec q_l, \vec q_r)\) with
\begin{equation}
  \label{eq:46}
  \tilde{\vec A} (\vec q_l, \vec q_r) = \vec A(\tilde{\vec q})\;, 
\end{equation}
then it is called a \emph{Roe mean value}
for \(\vec q_l, \vec q_r\).
For the Euler equations, Roe found the following consistent mean values:
\begin{equation}
  \label{eq:108}
  \begin{split}
    \tilde \rho & = \sqrt{\rho_l \rho_r}\;, \\
    \tilde u & = \frac{\sqrt{\rho_l} u_l + \sqrt{\rho_r}
      u_r}{\sqrt{\rho_l} + \sqrt{\rho_r}}\;, \\ 
    \tilde v & = \frac{\sqrt{\rho_l} v_l + \sqrt{\rho_r}
      v_r}{\sqrt{\rho_l} + \sqrt{\rho_r}}\;, \\ 
    \tilde w & = \frac{\sqrt{\rho_l} w_l + \sqrt{\rho_r}
      w_r}{\sqrt{\rho_l} + \sqrt{\rho_r}}\;, \\ 
    \tilde H & = \frac{\sqrt{\rho_l} H_l + \sqrt{\rho_r}
      H_r}{\sqrt{\rho_l} + \sqrt{\rho_r}}\;, \\ 
    \tilde c & = \sqrt{(\gamma - 1)\bigl(\tilde H -
      \tfrac{1}{2}\tilde{\vec v}^2\bigr)}
  \end{split}
\end{equation}
with~\(\tilde{\vec v}^2 = \tilde u^2 + \tilde v^2 + \tilde w^2\) in
the full three-dimensional case. For 2d, we simply omit the
values for the third velocity component~\(w\). 

A consequence of this is that the wave speeds in the Roe mean value
determine the numerical viscosity. Since in the sonic point this
turned out to be disadvantageous as upwinding in opposite direction
leads to a sonic glitch, Harten~\cite{harten-tvd} came up with the
following fix: He replaces the
absolute value of an eigenvalue \(\lambda\) of the Roe matrix
by\index{Harten fix}
\begin{equation}
  \label{eq:53}
  \phi(\lambda) =
  \begin{cases}
    \abs{\lambda} & \text{if}~\abs{\lambda} \geq \delta\;,\\
    (\lambda^2+\delta^2)/(2\delta) & \text{if}~\abs{\lambda} <
    \delta\;, 
  \end{cases}
\end{equation}
where \(\delta\) is a small parameter. The numerical viscosity
coefficient is here bounded below by \(\delta/2\).
Furthermore, we apply
the fix only to the acoustic waves. In the test cases performed below,
we found no visible difference in the numerical solutions.

\subsection{Achieving Mach number consistency}
\label{sec:achi-mach-numb}

Guillard and Viozat~\cite{viozat} show that for very small Mach
numbers the viscosity resulting from the wave-wise upwind as in the
standard Roe scheme leads to an incosistency: While in the low Mach
number limit pressure fluctuations scale with~\(\mathcal O(M^2)\), the
numerical scheme supports pressure fluctuations of
order~\(\mathcal O(M)\), where~\(M\) is the reference Mach number of
the flow field. Guillard and Viozat~\cite{viozat} identify the
numerical viscosity as the source of this inconsistency, whereas
Guillard and Murrone~\cite{murrone} use this as the starting point for
a preconditioner that modifies the numerical viscosities on the
different waves such that they are of the same order. The crucial
point is to reduce the viscosity on the acoustic waves to the same
order as on the advective waves, the shear and entropy waves.

Fleischmann et al.~\cite{fleischmann} recommence this idea and suggest
a simple low-dissipation modification for the acoustic waves to
achieve what in the following we call \emph{Mach number consistency}:
In the low Mach number limit, the viscosities on all waves are of the
same order. Note that this is a generalization of the original concept
since we allow also for higher viscosities as long as they are of the
same order for all waves. Keeping this terminology, Fleischmann et
al.\ provide two basic strategies to maintain Mach number consistency
in the Roe solver: a Mach number dependent upper bound for the
viscosity on the nonlinear, i.\,e.\ acoustic, waves or a lower bound
for the viscosity on the linear waves, i.\, e.\ entropy and shear
waves. With a fixed positive number~\(\phi\), the first approach,
which is the main point of~\cite{fleischmann}, leads to wave speeds
\begin{equation}
  \label{eqfl:1}
  \tilde\lambda_{1,4} = \tilde u \mp \min{\phi\abs{\tilde u},\tilde
    c}\;,\qquad \tilde\lambda_{2,3} = \tilde u\;,
\end{equation}
the second to
\begin{equation}
  \label{eqfl:3}
  \tilde\lambda_{1,4} = \tilde u \mp \tilde c\;,\qquad
  \tilde\lambda_{2,3} =  \sign (\tilde u) \max{\frac{\tilde
    c}{\phi},\abs{\tilde u}}\;. 
\end{equation}
It is easy to see that the approach resembles is a simplified model
for the incomplete Riemann solvers traditionally used as a means to
prevent the carbuncle. Only, for these the number~\(\phi\) is not
fixed but depends on some type of indicator. While in the first case,
for low Mach number flows, all eigenvalues of the viscosity matrix
are~\(\mathcal O(u)\), in the latter case they are of the
order~\(\mathcal O(c)\). The original purpose of Fleischmann et al.\
is to suggest the first approach as an alternative to the
second. Their idea is that its application would in some sense reduce
the transverse pressure fluctuations also in the vicinity of a shock
and, thus, help to prevent the carbuncle
phenomenon. In~\cite{FLEISCHMANN2020109762}, they suggest a modified
HLLC-scheme based on this approach.

One should at this point keep in mind that the second approach when
naively applied also may lead to unphysical
solutions~\cite{felix-hll-roe}, for purely low Mach number flows, the
numerical viscosity might even outweigh the physical viscosity in a
Navier-Stokes computation. So, it has to be handled with care.

On the other hand, the low dissipation approach might impair the
numerical stability of the resulting scheme. It is well known (and can
be found in standard text books, e.\,g.~\cite{toro,leveque,laney})
that in connection with the explicit Euler method as time Integration,
the viscosities obtained from the wave-wise application of standard
upwind are the lowest that would guarantee numerical stability. The
situation might improve if some higher order time integration, e.\,g.\
some Runge Kutta method, is employed.

The question now arises if the two approaches might be combined in
some way, i.\,e. imposing an upper bound to the acoustic speed and a
lower bound to the absolute value of the flow speed at the same time
while still maintaining Mach number consistency. For this purpose, we
suggest the following setting:
\begin{equation}
  \label{eqfl:4}
  \tilde\lambda_{1,4} = \tilde u \mp {\tilde c}^\beta
  \left(\min{\phi\abs{\tilde u},\tilde c}\right)^{1-\beta}\;, \qquad
  \tilde\lambda_{2,3} =  \left( \sign (\tilde u) \max{\frac{\tilde
    c}{\phi},\abs{\tilde u}}\right)^\beta {\tilde u}^{1-\beta}\;,
\end{equation}
which is kind of a weighted geometric mean between the less and the
more viscous approach that we call Mach number consistency. As
desired, the order of the numerical viscosities for Mach numbers
tending to zero is~\(\mathcal O(M^{2-\beta})\).

Since the low dissipation approach by Fleischmann et~al.\ aims at the
pressure perturbations perpendicular to the flow direction, while the
traditional approach solely focuses on the possibility of a shock
within the considered Riemann problem, it would be desirable
to adapt the weights~\(\beta\) and~\(1-\beta\) in a way that leads to
the low dissipation approach in weakly compressible Riemann problems
and the high dissipation approach in the presence of strong acoustic
waves. For this purpose, like in our earlier work on the
carbuncle~\cite{vietnam,kemm-lyon,habil-kemm,shallow-carbuncle,carb15},
we resort to the residual in the Rankine-Hugoniot condition relative
to the acoustic speed, i.\,e.\ in this case its Roe mean value. We
find for the residual
\begin{equation}
  \label{eqfl:113}
  \mathfrak r = \vec f(\vec q_r) - \vec f(\vec q_l) - \tilde u\, (\vec
  q_r - \vec q_l)
\end{equation}
and from that
\begin{equation}
  \label{eqfl:5}
  \frac{\mathfrak r}{\tilde c} = \lambda_4 \vec r_4 - \lambda_1 \vec r_1\;,
\end{equation}
which is the difference between the right and left running acoustic
waves. Instead of the more elaborate functions, which we have applied
in our earlier work, we just compute~\(\beta\) by
\begin{equation}
  \label{eqfl:6}
  \beta = \min{\log[10]{\max{\frac{\mathfrak r}{\tilde c},1}},1}\;,  
\end{equation}
which is easily coded in any programming language since the common
logarithm~\(\log[10]{}\) is usually available as an intrinsic function
in the compilers as well as the natural logarithm. Since the residual
vanishes for shear and entropy waves and is rather large for strong
shocks, this choice of the weight~\(\beta\) ensures that in the low
Mach number regime we are left with the low dissipation approach,
while a shock in the Riemann problem enforces the high dissipation
Mach number consistent scheme.

Since the powers~\((\cdot)^\beta\) and~\((\cdot)^{1-\beta}\) in
equation~\eqref{eq:4} are computationally expensive on most machines,
in our numerical investigation, we also include an approximation of
the weighted geometric mean: the corresponding weighted arithmetic
mean
\begin{equation}
  \label{eq:8}
  \tilde\lambda_{1,4} = \tilde u \mp \beta \tilde c + (1-\beta)
  \min{\phi\abs{\tilde u},\tilde c}\;, \qquad
 \tilde\lambda_{2,3} = \beta \sign (\tilde u) \max{\frac{\tilde
    c}{\phi},\abs{\tilde u}} + (1-\beta)\tilde u\;.
\end{equation}
Although this would not guarantee full Mach number consistency, it
might still be a reasonable replacement at lower numerical cost.

\section{Numerical investigation}
\label{sec:numer-invest}

In order to assess the behaviour of the simple prototypical Mach
consistent Roe solvers discussed above, we perform a series of
numerical tests. The test cases are chosen such that we can shed light
on both the resolution of the resulting scheme and the
stability of discrete shock profiles. Furthermore, we hope to find
information on the possible loss in terms of numerical stability that
might result from lowering the viscosity on the acoustic
waves.

\subsection{Overview of the test cases}
\label{sec:overview-test-cases}

Some of the tests can be performed in one space dimension
while others are genuinely two-dimensional. For the sake of
simplicity, we restrict ourselves to two dimensions and, thus, do not
perform three-dimensional tests.

\subsubsection{Steady shear wave}
\label{sec:steady-shear-wave}

In order to investigate the resolution of the schemes, it is useful to
consider a simple steady shear wave. A complete Riemann solver is
expected to resolve steady shear waves exactly. In our tests, we
modify the initial data by adding randomized numerical noise. While
for the solver with the viscosity on shear and entropy waves bounded
below, we expect the wave to be smeared out to some extent and for the
solver with the viscosity on acoustic waves bounded to be exactly
resolved, the main question is how the presence of the perturbations
will affect the behaviour of the blended schemes. Is it sufficient to
employ an indicator function in order to obtain results of the same
quality as with the pure low viscosity scheme? The unperturbed initial
data would be~\(\rho =1,\ p=1,\ u=0,\ v_{r/l} = \pm 1\). The amplitude
of the artificial noise is~\(10^{-6}\). We will perform the test in
two versions: purely one-dimensional and on the other hand
two-dimensional where the shear wave is located on a grid line. This
way, the 2d-version can be considered as a set of---except for the
perturbations---identical 1d problems which somehow interact.
This consideration will also be adopted when plotting the results of
above tests.

\subsubsection{Colliding flow}
\label{sec:colliding-flow}

This test~\cite[Section~7.7]{astro-leveque} resembles a simplified
model for the starting process of the blunt body test when using the
inflow state as initial data in the complete computational
domain. This is best understood when considering the flow before the
blunt body along the symmetry line. Since the flow is aligned with
that symmetry line, and due to the switch of the sign of the flow
velocity for wall boundary conditions, it behaves essentially like the
left half of a colliding flow in 1d. To turn it into a 2d-test, the
flow is, again, equipped with an additional space direction, in which
everything is expected to be constant. In order to trigger the
carbuncle, the initial state is superimposed with noise that is
generated randomly and has a small amplitude.

For our numerical test, in the initial state, density and pressure are
set to~\(\rho=1,\ p=1\).  The normal velocity is set
to~\(u_\text{left/right} = \pm 20\), the transverse velocity component
to~\(v=0\).  To trigger the carbuncle, we superimpose artificial
numerical noise of amplitude \(10^{-6}\) onto the primitive variables
instead of disturbing it in just one point as was done originally by
LeVeque~\cite[Section~7.7]{astro-leveque}. The 2d computations are
done on the rectangle~\([0,60]\times[0,30]\) discretized with
\(60\times 30\) grid cells. For the 1d computations, we just drop the
second space dimension.

Since in 2d, the problem is a simple 2d-extension of the one-dimensional
problem, the results are presented in scatter-type plots: we slice the
grid in \(x\)-direction
along the cell faces and plot the density for all slices at once.

\subsubsection{Uniform flow}
\label{sec:uniform-flow}

In its nature, this problem is similar to the two examples above: it
is a 1d problem, which we artificially extend to two space
dimensions. The basic flow is a uniform flow in \(x\)-direction but
superimposed with artificial random numerical noise. We consider two
examples, one with Mach number~\(M=20\) and one with Mach
number~\(M=1/20\). Thus, we test the behaviour for highly supersonic as
well as for weakly compressible flows. The unperturbed initial values
for density and flow velocity (in \(x\)-direction)
are~\(\rho = 1,\ u=1\). The pressure is adjusted to the desired Mach
number.  The amplitude of the artificial noise is~\(10^{-6}\).

\subsubsection{Steady shock}
\label{sec:steady-shock}

While the colliding flow test models the starting process of the blunt
body flow, the steady shock test, introduced by Dumbser
et~al.~\cite{michael-carbuncle}, features a simplified model for the
converged shock in the blunt body flow. Following Dumbser et~al., we
set in the upstream region~\(\rho =1,\ u=1\).  The upstream Mach
number is set to~\(M=20\), the transverse velocity component
to~\(v=0\), which again results in the 2d extension of a 1d problem.
The shock is located directly on a cell face.  To trigger the
instability of the discrete shock profile, we add artificial numerical
noise of amplitude \(10^{-6}\) to the primitive variables in the
initial state. The computations are done on~\([0,100]\times[0,40]\)
discretized with~\(100\times 40\) grid cells.

Again for the presentation of the results, we also employ scatter-type
plots as described for the colliding flow
problem~\ref{sec:colliding-flow}.

\subsubsection{Quirk test}
\label{sec:quirk-test}

Quirk~\cite{quirk} introduced a test problem which is known as
Quirk test. Contrary to the steady shear, the colliding flow, and the
steady shock, it is not a
one-dimensional Riemann problem, but consists of a shock running down
a duct. The shock is caused by Dirichlet-type boundary conditions on
the left boundary with~\(\rho=5.26829268\), \( u=4.86111111\),
\(p=29.88095238\), while the flow field is initialized with~\(\rho = 1,\
u=v=0,\ p=1/\gamma\). Originally, a disturbance of the middle grid
line was used to trigger the instability~\cite{quirk}. Because the
computations are done with a Cartesian code, we instead use numerical
noise in the same manner as for the steady shock and the colliding
flow problem. The only difference lies in the amplitude of the
perturbation, here \(10^{-3}\). The computations are done
on~\([0,1600]\times[0,20]\) discretized with \(1600\times 20\) grid
cells. 
Again we use scatter-type
plots (as described above) to present the results.

\subsubsection{Elling test}
\label{sec:elling-test}

Elling~\cite{elling_carbuncle_2009} investigates the influence of the
supersonic upstream region on the shock profile. He models the
interaction of a vortex filament with a strong shock. For this
purpose, he starts with a steady shock. In the region upstream of the
shock, he picks the middle slice of the computational grid and
artificially sets the velocity to zero. For more details, we refer to
our previous work~\cite{carb15}.
The
initial condition is a modified version of the
steady shock test, cf.\ Section~\ref{sec:steady-shock}. The
region to the right of the shock remains unchanged. In the supersonic
inflow region, only the middle \(x\)-slice is changed. Here, the velocity is
set to zero. This is done to model a vortex layer hitting the shock
front, which in turn is a prototype for shock-boundary-layer
interaction. 
%

\subsubsection{Kelvin-Helmholtz instability}
\label{sec:kelv-helmh-inst}

For the Kelvin-Helmholtz instability, we start with a flow which
consists of three parts. In a region close to the middle line, it is
uniformly directed towards the left with~\(u=-\sqrt{\gamma}/2\). Above
and below, it is directed in the opposite direction
with~\(u=\sqrt{\gamma}/2\). Pressure and density are constant
everywhere with~\(p=1\) and~\(\rho=1\). To trigger the instability, we
slightly disturb the \(y\)-component of the flow velocity. The
perturbation is sinusoidal in~\(x\)-direction with an amplitude
of~\(1/100\) and constant in~\(y\)-direction. 
However, as we cannot expect for a first order scheme that the
instability will evolve on our given~\(100\times 100\) grid, we resort
to a second order scheme with direction-wise minmod on primitive
variables.
Since this example is dominated by entropy and shear waves, it is a
good test for the overall resolution of the scheme. Furthermore, it
features a low Mach number flow, which allows us to evaluate the
robustness of the scheme and how much fragility might arise from the
upper bound for the viscosity on the acoustic waves.

\subsubsection{Double Mach Reflection (DMR)}
\label{sec:double-mach-refl}

The Double Mach Reflection is a self similar solution that evolves
from a shock running up a ramp. It was introduced by Woodward and
Colella~\cite{WC} as a benchmark for Euler codes. In contrast to the
single Mach reflection, the double Mach reflection features two triple
points and, thus, is a challenging problem for gas dynamics
codes. Here, only a brief overview is given since this problem has
already been discussed in more detail in or earlier
work~\cite{dmr}. As suggested therein, we double the size of the
computational domain in the \(y\)-direction and also follow the other
suggestions made in that study. The results are displayed as contour
plots of the transverse momentum~\(\rho v\), which gives the best
insight into the quality of the scheme. This is especially true for
the resolution of the secondary slip line.

\subsection{Numerical results and discussion}
\label{sec:numer-results-disc}

For the numerical results, we employed Euler2d, a simple 2d-Cartesian code
developed for the test of Riemann solvers in the Group of Claus-Dieter
Munz at Stuttgart University. The code implements standard finite
volumes.  For higher order, direction-wise geometric limiting with
minmod on primitive variables is used.

\subsubsection{One-dimensional tests}
\label{sec:one-dimens-tests}

In this section, we consider two tests: a perturbed steady shear wave
and a colliding flow. While the purpose of the shear wave test is
obvious, the purpose of the latter is to check for issues, that the
lowered viscosity on the acoustic waves might cause.
\begin{figure}
  \centering
  \includegraphics[width=.9\linewidth]{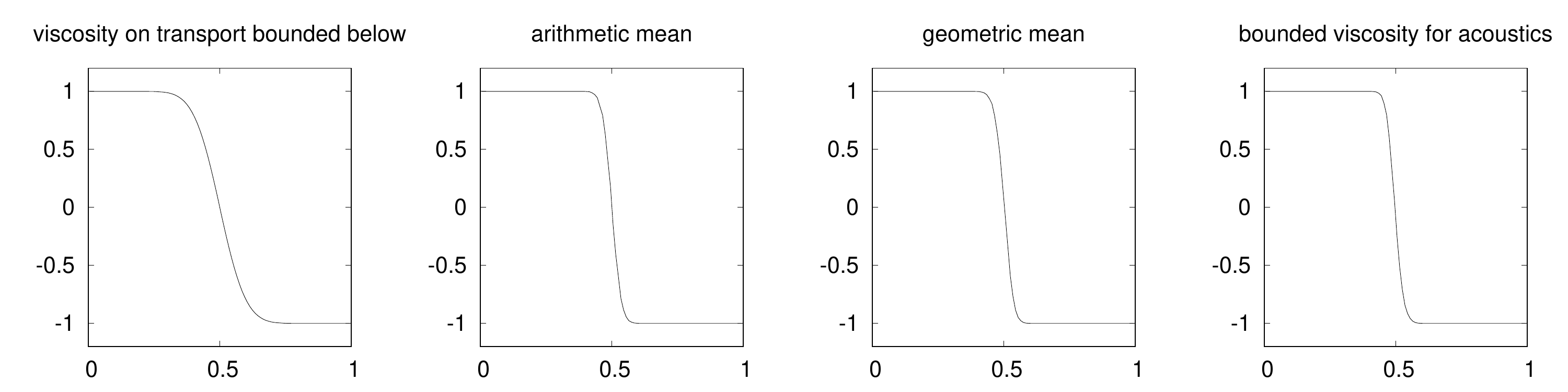}
  \caption{Transverse velocity for steady shear wave with randomized
    noise at time~\(t=2.5\)}
  \label{fig:shear1d}
\end{figure}
The results for the \textbf{perturbed steady shear wave} in
Figure~\ref{fig:shear1d} clearly show that the indicator
function~\(\beta\) nicely detects the absence of shocks and, thus,
avoids unnecessary viscosity on the shear wave. Except for the high
viscosity approach, all versions show the same result.

\begin{figure}
  \centering
  \includegraphics[width=.9\linewidth]{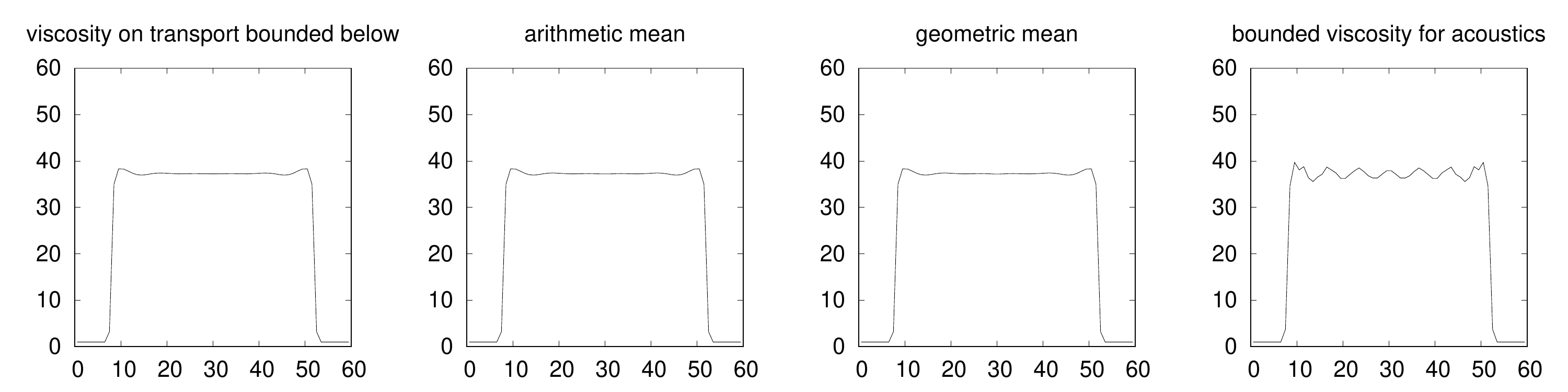}
  \caption{Colliding flow at time~\(t=30\)} 
  \label{fig:colliding1d}
\end{figure}
In Figure~\ref{fig:colliding1d}, for the \textbf{colliding flow} test,
we find the opposite situation: here, all versions except for the
original Fleischmann solver yield the same results. If the viscosity
on the acoustic waves is bounded regardless of the presence or absence
of strong acoustic waves, we suffer from a strong increase of
post-shock oscillations.


\subsubsection{2d counterparts of 1d tests and other quasi one domensional
  2d tests}
\label{sec:2d-counterparts-1d}

Since these tests are simple 2d-extensions of one-dimensional problems,
the results are presented in scatter-type plots: we slice the grid in
\(x\)-direction along the cell faces and plot the density or, in the
case of the perturbed shear wave, the transverse velocity for all
slices at once.

\begin{figure}
  \centering
  \includegraphics[width=.9\linewidth]{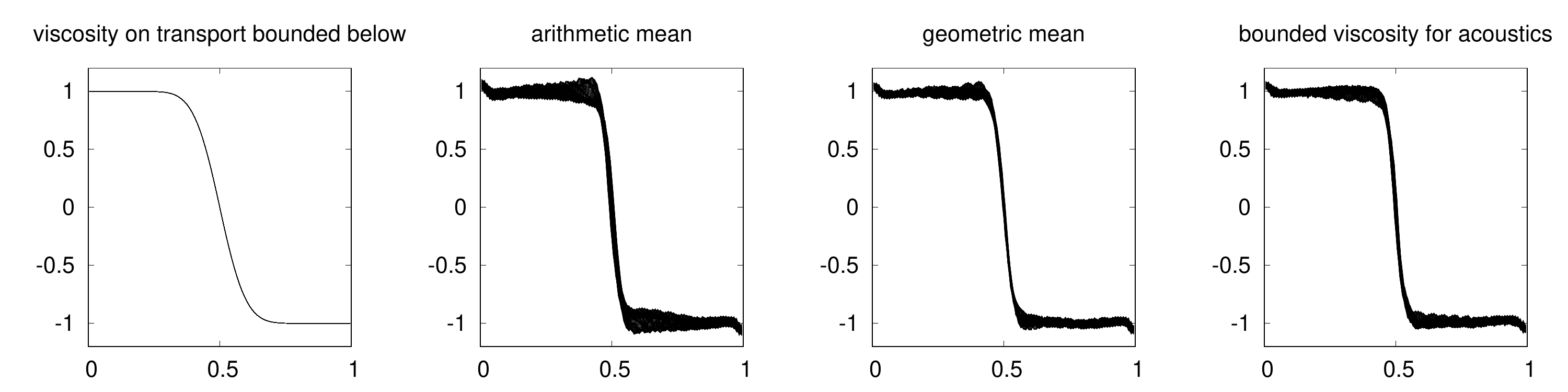}
  \caption{Steady shear wave with randomized noise at time~\(t=2.5\)} 
  \label{fig:shear2d}
\end{figure}

As can be seen in Figure~\ref{fig:shear2d}, the 2d \textbf{perturbed
  steady shear wave} reveals that the strict application of a lower
bound for the transport viscosity at the one hand leads to a strong
smearing effect on the shear wave but on the other hand prevents
perturbations in \(y\)-direction. Furthermore, we indeed observe a
slight difference between the application of the arithmetic and the
geometric mean. As we will see, this is typical for situations without
strong acoustic waves.

\begin{figure}
  \centering
  \includegraphics[width=.9\linewidth]{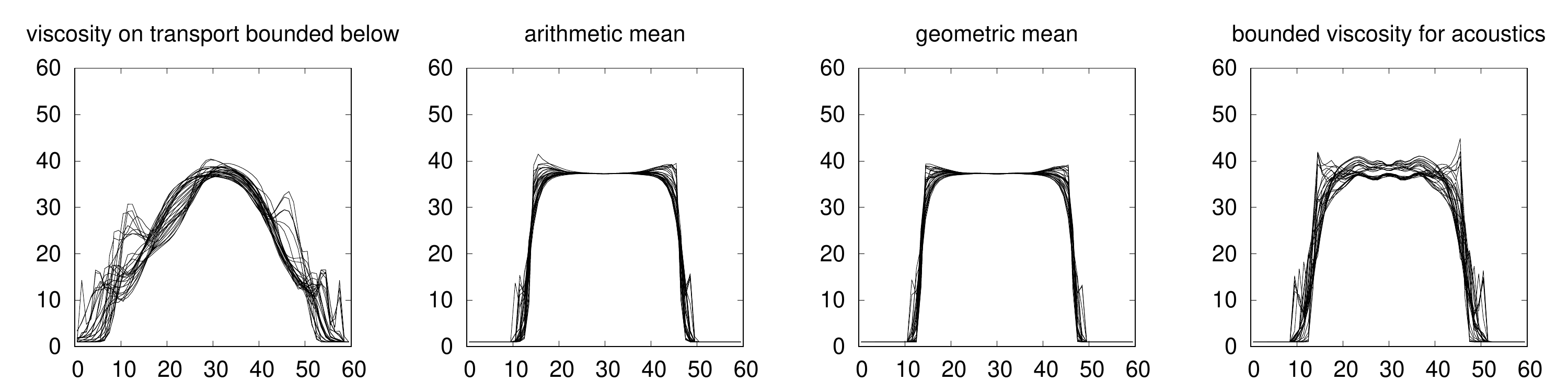}
  \caption{Colliding flow at time~\(t=30\)} 
  \label{fig:colliding2d}
\end{figure}

For the 2d \textbf{colliding flow}, as depicted in
Figure~\ref{fig:colliding2d}, we get quite reasonable, but not
yet perfect, results with the blended schemes. An interesting point is
that the low dissipation version yields better results than the high
viscosity variant. This is surprising, especially with respect to the
purely one-dimensional results in Section~\ref{sec:one-dimens-tests}. 

\begin{figure}
  \centering
  \includegraphics[width=.9\linewidth]{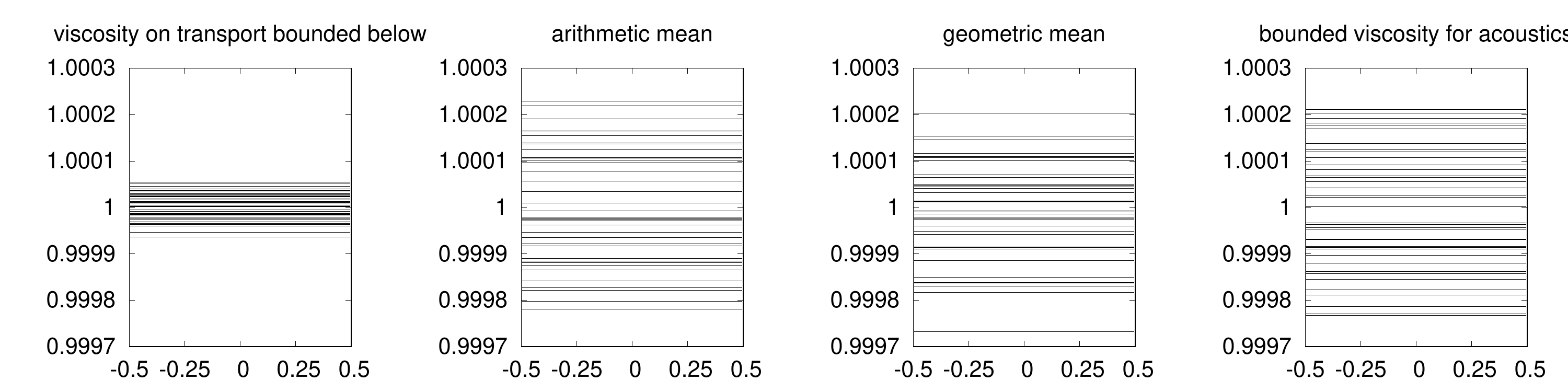}
  \caption{Supersonic uniform flow with randomized noise at time~\(t=5\)} 
  \label{fig:unisuper}
\end{figure}
\begin{figure}
  \centering
  \includegraphics[width=.9\linewidth]{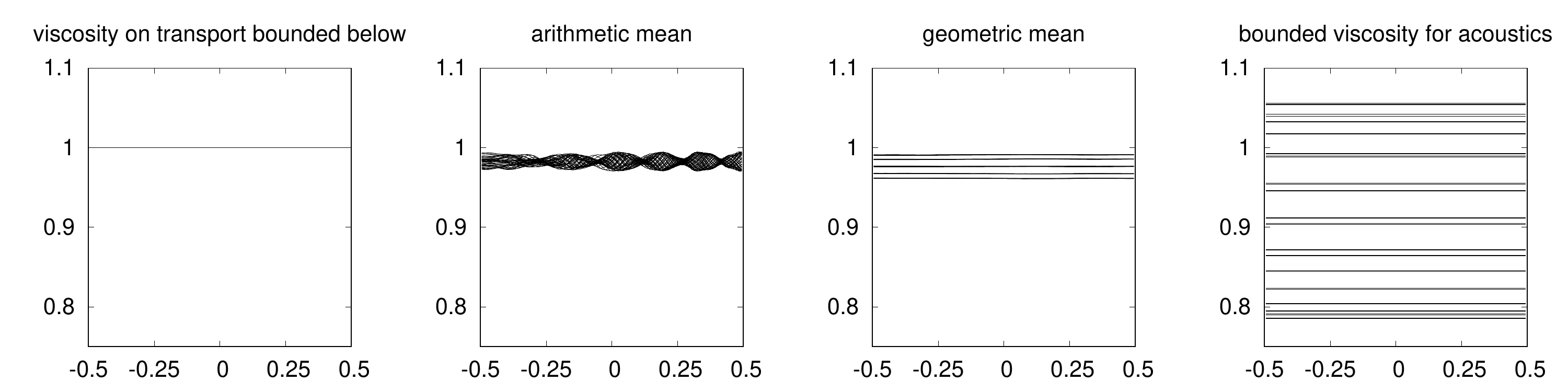}
  \caption{Subsonic uniform flow with randomized noise at time~\(t=5\)} 
  \label{fig:unilowmach}
\end{figure}

For the \textbf{uniform flow}, we present results for both
the~\(M=20\) and the~\(M=1/20\) flow in Figures~\ref{fig:unisuper}
and~\ref{fig:unilowmach}, respectively. While for the former the
differences between the blended and the low dissipation schemes are
rather small, for the latter they are significant. Also, the
differences between the geometric and the arithmetic mean are clearly
visible. Although the amplitude of the perturbations is about the
same, the version with the geometric mean shows only variations in the
transverse direction while the variant with the arithmetic mean
features a more complex flow pattern.

\begin{figure}
  \centering
  \includegraphics[width=.9\linewidth]{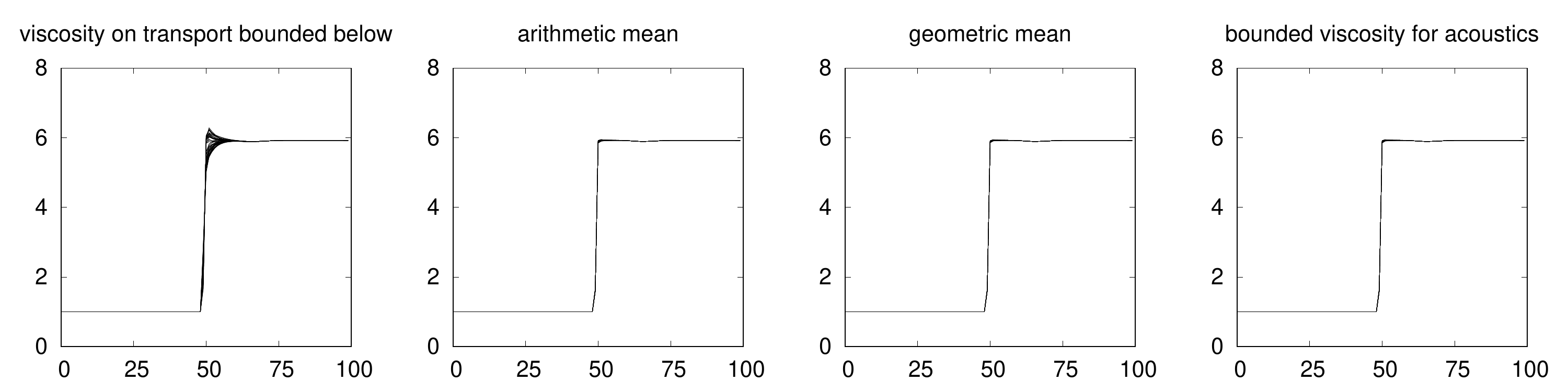}
  \caption{Steady shock test at time~\(t=100\)} 
  \label{fig:steady}
\end{figure}

The test case with the \textbf{steady shock} is well resolved by all
solvers, cf.\ Figure~\ref{fig:steady}. With the high viscosity method,
which is intended to serve as a simplified model for the traditional
carbuncle cures, the results are quite similar to those of the
carbuncle cures which rely on an indicator computed from the Riemann
problem itself, e.\,g.\ our own HLLEMCC solver~\cite{kemm-lyon}. Even
the blended and the low dissipation schemes yield perfect
results. This confirms that there indeed is, as stated by Fleischmann
et al.~\cite{fleischmann}, a significant influence of the post shock
region on the stability of the discrete shock profile. Furthermore, it
seems even better to opt for the low dissipation Mach consistent
numerical flux than for the traditional high dissipation approach.

\begin{figure}
  \centering
  \includegraphics[width=.9\linewidth]{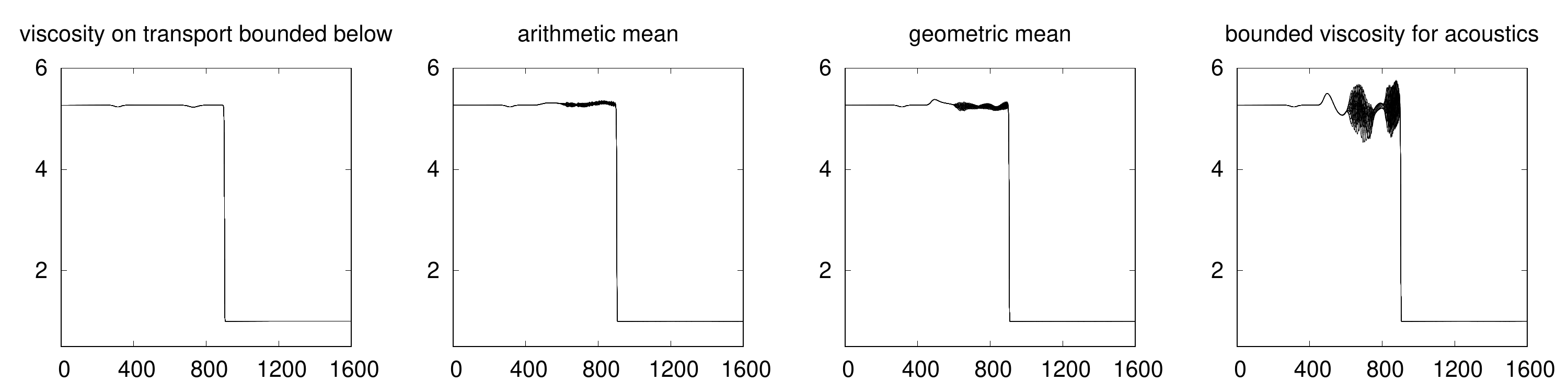}
  \caption{Quirk test at time~\(t=150\)} 
  \label{fig:}
\end{figure}
As can be seen in Figure~\ref{fig:}, the results for the \textbf{Quirk
  test} (cf.\ Section~\ref{sec:quirk-test}) show some similarity with
the well known shock entropy wave interaction problem as suggested by
Shu and Osher test~\cite{shu-osher}, where a moving shock interacts
with density fluctuations of small amplitude. While for the version
with the transport viscosity bounded below the results look perfect,
the effect is clearly visible for the strict low dissipation Roe
solver. For the blended schemes, it is still visible but with a much
smaller amplitude. This indicates that the lowered viscosity on the
acoustics in the resting fluid right of the shock causes slight
numerical instabilities which lead to an increase of the density
perturbations introduced by the artificial noise in the initial state.

Furthermore,  as for the low Mach number uniform flow in
Figure~\ref{fig:unilowmach}, we observe a difference between the
blended versions with geometric and arithmetic mean. This might be due
to the fact that the perturbed resting fluid resembles the situation
of the Mach number uniform flow closely. In fact, it is a flow field
with extremely low Mach number.


\subsubsection{Genuinely two-dimensional test cases}
\label{sec:genu-two-dimens}

Here, we discuss some problems that are genuinely two-dimensional,
which means that they have no one-dimensional counterpart. They allow
us to study how the results we have gotten so far adapt to more
real-life situations with complex flow structures. Since the problems
are genuinely 2d, we also have to resort to genuinely 2d-types of
plots. In most cases, we show contour plots, while for one example a
surface plot is also provided since the contour plot alone would be
somehow misleading.
\begin{figure}
  \centering
  \includegraphics[width=.9\linewidth]{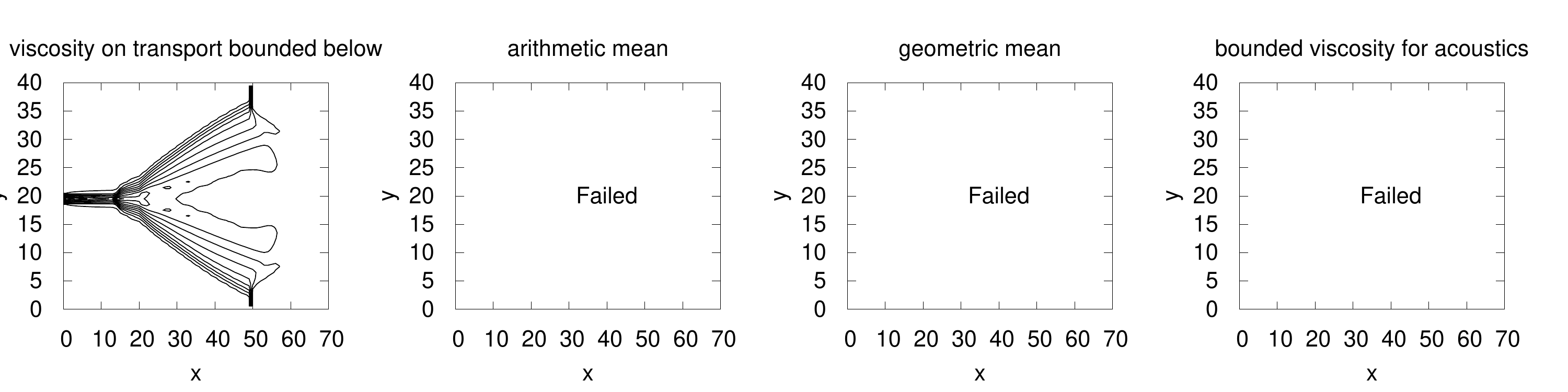}
  \caption{Density for Elling test at time~\(t=100\); 1st order computation} 
  \label{fig:elord1}
\end{figure}

\begin{figure}
  \centering
  \includegraphics[width=.9\linewidth]{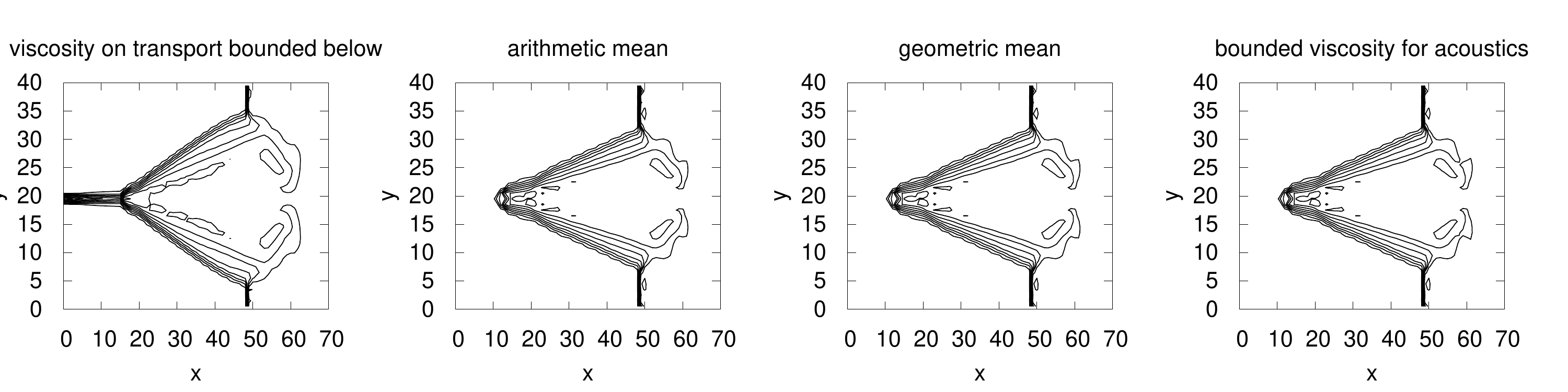}
  \caption{Entropy for Elling test at time~\(t=100\); 2nd order computation} 
  \label{fig:elord2}
\end{figure}

In the \textbf{Elling test} (Figure~\ref{fig:elord1}), we can see the
stability issues introduced by the low dissipation on the acoustic
waves. All first order computations except for the high dissipation
method abort after some time, which is however fixed by the increase
of the order of the scheme, as can be seen in
Figure~\ref{fig:elord2}. This also might explain why in their studies
Fleischmann et al.~\cite{fleischmann,FLEISCHMANN2020109762} did not
observe any stability issues: they employed fifth order in space and
third order in time. For the SSP-Runge-Kutta, they use for the time
integration, no details are given by Fleischmann
et~al.~\cite{fleischmann}. But there is a good chance that also the
stability region is shaped in a way that provides additional
stability. As a comparison to the results in~\cite{carb15} shows,
while the result of the high dissipation closely resembles our old
carbuncle cure, HLLEMCC~\cite{kemm-lyon}, the other versions resemble
the results obtained with the Osher solver.

\begin{figure}
  \centering
  \includegraphics[width=.9\linewidth]{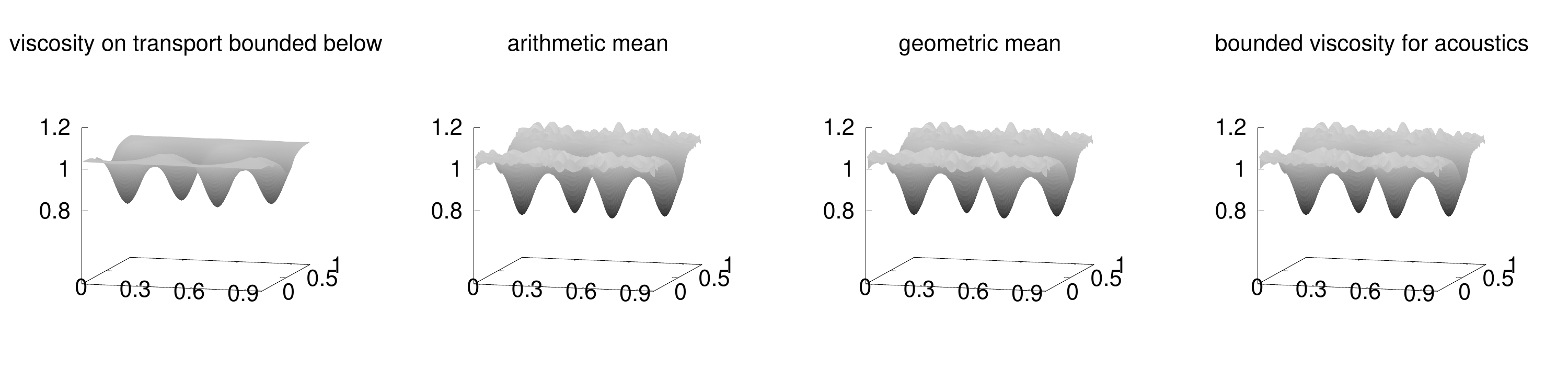}
  \caption{Entropy for Kelvin-Helmholtz instability at time~\(t=4\),
    surface plot} 
  \label{fig:KHIpmd3}
\end{figure}

\begin{figure}
  \centering
  \includegraphics[width=.9\linewidth]{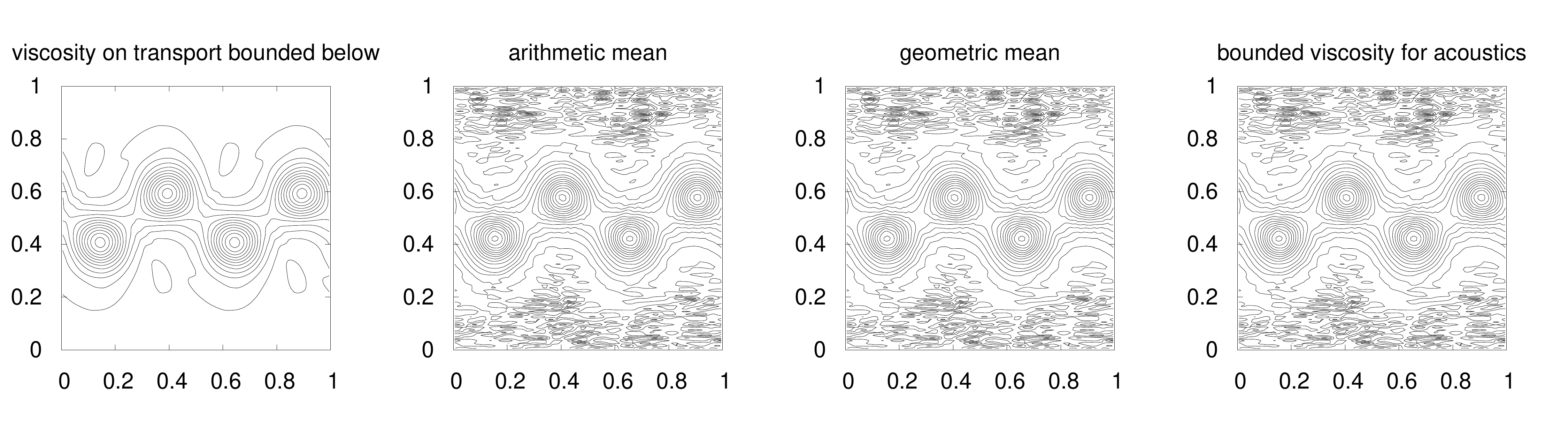}
  \caption{Density for Kelvin-Helmholtz instability at time~\(t=4\),
    contour plot} 
  \label{fig:KHIcontour}
\end{figure}

The results for the \textbf{Kelvin-Helmholtz} instability (cf.\
Section~\ref{sec:kelv-helmh-inst}) show a better resolution of the
instability itself when the transport viscosity is not bounded below
in weakly compressible flows (Figure~\ref{fig:KHIpmd3}). But on the
other hand, and much more prominent, at least when the contours are
plotted as in Figure~\ref{fig:KHIcontour}, some instabilities indicate
that the naive approach for lowering the viscosity on acoustic waves
as studied in this paper is not sufficient. Like for the low Mach
uniform flow and the Quirk test as well as the Elling test, one would
ask for a more elaborate way of adjusting the acoustic viscosity,
maybe by adapting the parameter~\(\phi\) in some way.

\begin{figure}
  \centering
  \includegraphics[width=.9\linewidth]{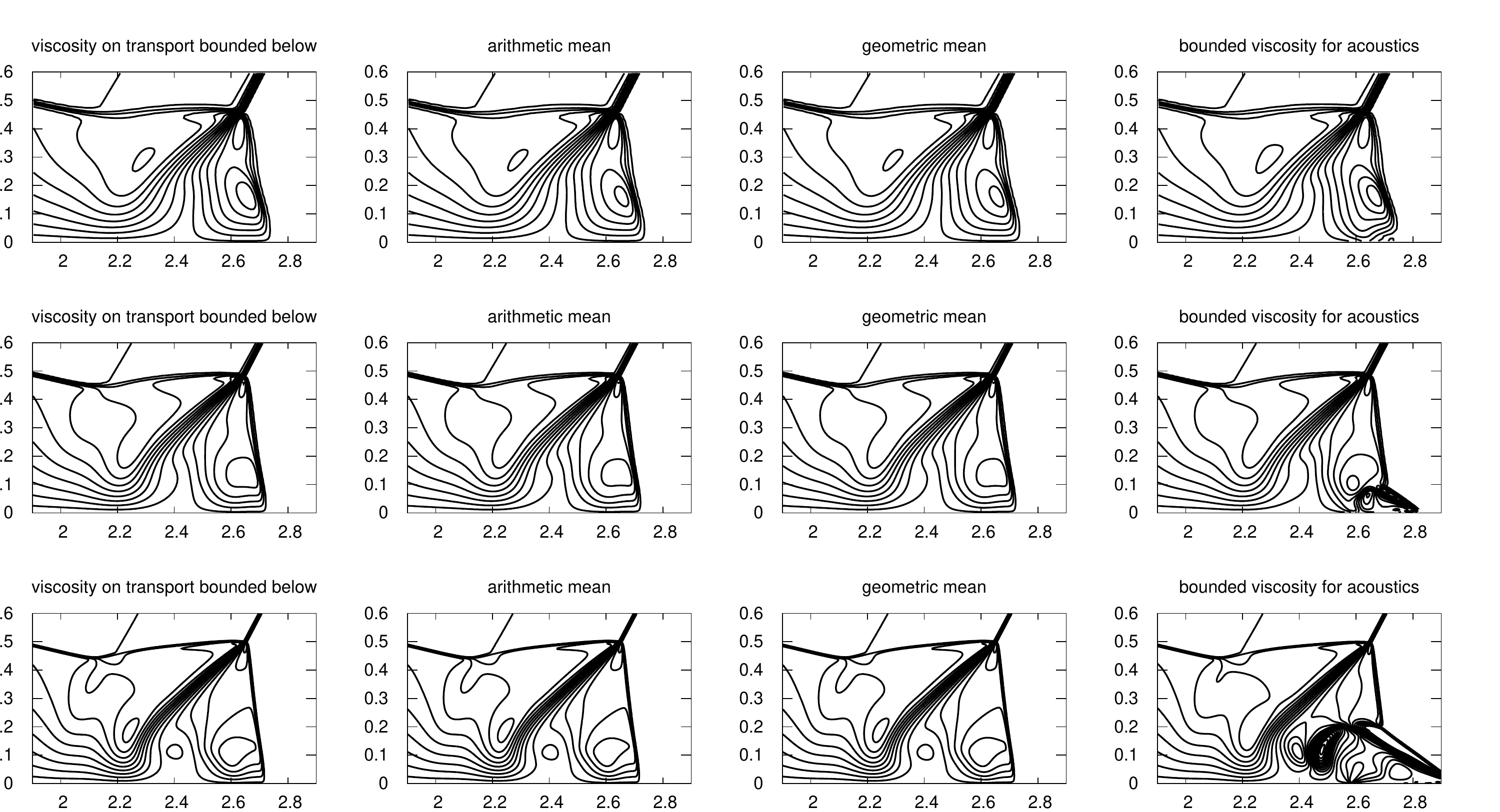}
  \caption{Vertical momentum for Double Mach Reflection problem; 1st
    order computations with~\(\Delta x = \Delta y =
    1/120\),~\(\Delta x = \Delta y = 1/240\),
    and~\(\Delta x = \Delta y = 1/480\) (from top to bottom) }
  \label{fig:dmrord1}
\end{figure}

\begin{figure}
  \centering
  \includegraphics[width=.9\linewidth]{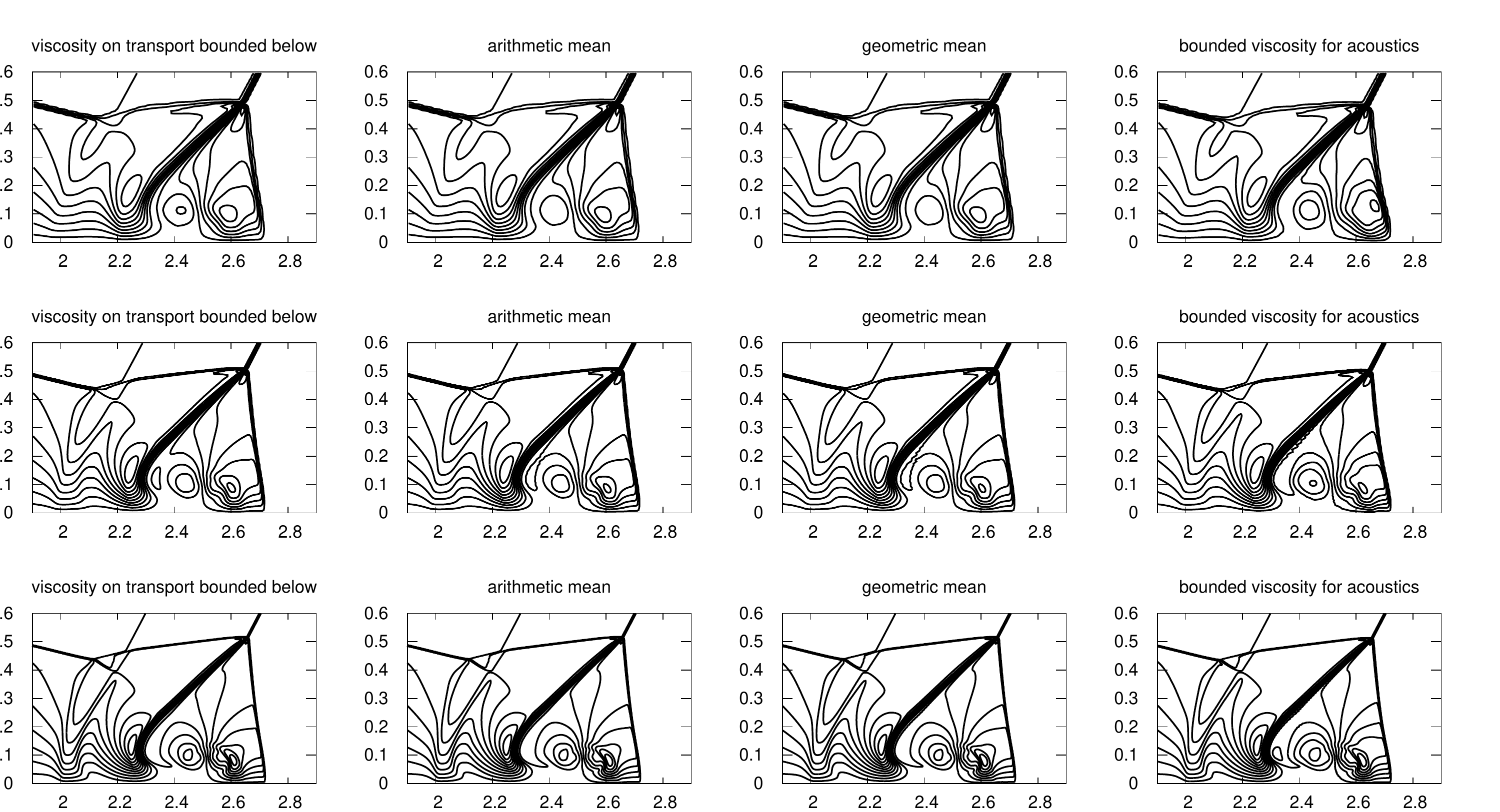}
  \caption{Vertical momentum for Double Mach Reflection problem; 2nd
    order computations with~\(\Delta x = \Delta y =
    1/120\),~\(\Delta x = \Delta y = 1/240\),
    and~\(\Delta x = \Delta y = 1/480\) (from top to bottom)} 
  \label{fig:dmrord2}
\end{figure}

Finally, we show results for the \textbf{double Mach reflection} as
described in Section~\ref{sec:double-mach-refl}. While
Figure~\ref{fig:dmrord1} displays the results of the first order
computation, Figure~\ref{fig:dmrord2} provides the results obtained
with second order in space and time. We performed computations with
a grid spacing of~\(\Delta x = \Delta y =
1/120\),~\(\Delta x = \Delta y = 1/240\),
and~\(\Delta x = \Delta y = 1/480\), where the latter is the standard
resolution for the DMR test as introduced by Woodward and
Colella~\cite{WC}. The most prominent feature in the results is the
kinked Mach stem originating from a carbuncle in the computations with
the pure low dissipation approach. This indicates that, in order to
safely prevent a carbuncle, it is still necessary to increase the
viscosity on shear and entropy waves in the vicinity of strong
shocks. Fighting one possible source of the carbuncle seems not 
sufficient. In fact, it seems to be crucial to address all possible
sources at the same time in order to get optimal results.

For the second order computations the situation is much better, as no
carbuncle can be found. Furthermore, the differences between the
schemes are much smaller, and the secondary slip line can be seen in
all of them, especially with~\(\Delta x = \Delta y =
1/480\). According to our earlier findings in~\cite[\S\,5.2]{carb15},
the higher order already stabilizes the shock position. Although
minmod is rather weak in this respect, for the DMR, it seems sufficient to provide the necessary additional stabilization.

\section{Conclusions and possible directions for further research}
\label{sec:concl-poss-direct}

We started with the work of Fleischmann et~al.~\cite{fleischmann} and
their simplified models for low Mach number corrections~\eqref{eqfl:1}
and carbuncle preventing incomplete Riemann solvers~\eqref{eqfl:3} and
combined them via the blending~\eqref{eqfl:4}
or its computationally cheaper approximation~\eqref{eq:8} in
connection with the indicator~\eqref{eqfl:6} to
investigate the influence of Mach number consistency on the carbuncle
phenomenon. The parameter~\(\beta\) was introduced to provide a smooth
transition between the low dissipation scheme for weakly compressible
and the high dissipation method for fully compressible Riemann
problems. It should be kept in mind that in high-speed grid aligned
flows both cases occur at the same time, depending on the space
direction.

The first results look rather promising. For many test cases, at least
when the order of the scheme is increased, the schemes perform
well. For the steady shock problem, the perturbations along the shock
front were much smaller than with traditional carbuncle cures based
solely on the use of a higher viscosity for shear and entropy waves
near shocks. In some test cases like the Elling test, where all except
the high dissipation scheme failed in the 1st order computations, as
well as the Quirk test (mainly the pure low dissipation version) and
the Kelvin-Helmholtz instability, perturbations caused by the
violation of the stability conditions as a result of lowering the
viscosity below the Roe values can be seen. Addressing this issue
seems to be the most important contribution expected from future
research.

Since there is ongoing and vivid research towards Riemann solvers for
low Mach number flows,
e.\,g.~\cite{felix-low-mach-fix,philipp-l2roe,pascal-et-al-low-mach,Barsukow2017,GUILLARD2017203,CHEN20183737,all-mach-hllc,FLEISCHMANN2020109762},
there is also hope that in the near future a way will be found to
construct a robust all Mach number Riemann solver that prevents
unphysical carbuncles based on the above considerations. Although Xie
et~al.~\cite{all-mach-hllc} do not explicitly control Mach number
consistency, their solver can already be considered a step in this
direction.

In some cases, the experimental solver used in this study might be
sufficient: High order SSP/monotone time integration schemes, and many
high order Runge-Kutta schemes in general, tend to provide some sort
of reserve in terms of stability, something the explicit Euler scheme
lacks for all modes of the semi-discrete system obtained by first
order standard upwind. As the original results by Fleischmann
et~al.~\cite{fleischmann} and the results for their modified
HLLC-solver~\cite{FLEISCHMANN2020109762} show, for most computations,
this seems to be sufficient even when the low dissipation
version~\eqref{eqfl:1} is applied in its pure form. Thus, we conclude
that the blended schemes can also be safely used in many very high
order gas dynamics codes.

\bibliographystyle{amsplain} \bibliography{carb-fl}

\end{document}